\documentclass[a4paper,11pt,twoside]{article}
\usepackage{pst-fill,pst-grad}
\usepackage{textcomp}
\usepackage[english]{babel}
\usepackage{graphicx}
\usepackage{amsmath}
\usepackage{float}
\usepackage[matrix,arrow,curve]{xy}
\usepackage{pstricks}
\usepackage{amsmath,amsfonts,verbatim,afterpage,theorem,euscript,mathrsfs,amssymb}
\usepackage{amsfonts}
\usepackage{amssymb}
\usepackage{array}
\usepackage{dsfont}
\usepackage{hyperref}
\newtheorem{Theoreme}{Theorem}

\newtheorem{Proposition}{Proposition}[section]
\newtheorem{Lemme}{Lemma}[section]

\newtheorem{Remarque}{Remark}[section]

\numberwithin{equation}{section}
\def\vu{\vec{u}}
\def\vv{\vec{v}}
\def\vf{\vec{f}}
\def\vg{\vec{g}}

\def\vn{\vec{\nabla}}

\def\Rt{\mathbb{R}^3}

\def\mPl{\mathcal{P}^{\text{log}}}

\title{\bf  Lebesgue spaces with variable exponent: some applications to the Navier-Stokes equations} 
\usepackage{authblk}
\author{Diego Chamorro\footnote{\emph{diego.chamorro@univ-evry.fr} (corresponding author)} }
\author{Gast\'on Vergara-Hermosilla\footnote{\emph{gaston.vergarahermosilla@univ-evry.fr}} }
\affil{\footnotesize LaMME, Univ. Evry, CNRS, Universit\'e Paris-Saclay, 91025, Evry, France.}
\begin{document}
\maketitle
\begin{scriptsize}
\abstract{In this article we study some problems related to the incompressible 3D Navier-Stokes equations from the point of view of Lebesgue spaces of variable exponent. These functional spaces present some particularities that make them quite different from the usual Lebesgue spaces: indeed, some of the most classical tools in analysis are not available in this framework. We will give here some ideas to overcome some of the difficulties that arise in this context in order to obtain different results related to the existence of mild solutions for this evolution problem. }\\[3mm]
\end{scriptsize}
%\tableofcontents
%%%%%%%%%%%%%%%%%%%%%%%%%%%%%%%%%%%%%%%%%%%%%%%%%%%
\section{Introduction}
We consider here the classical incompressible Navier-Stokes equations defined in the whole space $\mathbb{R}^3$:
\begin{equation}\label{NS_Intro}
\begin{cases}
\partial_t\vu=\Delta \vu-(\vu \cdot \vn) \vu-\vn P + \vec{f}, \qquad div(\vu)=0,\\[5pt]
\vu(0,x)=\vu_0(x),\quad div(\vu_0)=0, \qquad x\in \mathbb{R}^3,
\end{cases}
\end{equation}
where $\vu:[0, +\infty[\times \mathbb{R}^3 \longrightarrow \mathbb{R}^3$ is the velocity field, $P:[0, +\infty[\times\mathbb{R}^3 \longrightarrow  \mathbb{R}$ denotes the pressure and $\vu_0:\mathbb{R}^3 \longrightarrow \mathbb{R}^3$, $\vf:[0, +\infty[\times \mathbb{R}^3 \longrightarrow \mathbb{R}^3$ are a given initial data and a given exterior force, respectively.\\

From the point of view of the existence of solutions we have at our disposal (at least) two main theories: \emph{mild} solutions which are local in time (for any generic initial data $\vu_0$) and \emph{weak} (Leray) solutions which satisfy an energy inequality and are global in time.\\

In this article we are concerned by the existence of \emph{mild} solutions which are obtained with the help of a fixed-point theorem. Of course, in this theory the choice of a \emph{good} functional setting is crucial: indeed, from the seminal work of Fujita-Kato \cite{fujita1964navier} (where classical Lebesgue spaces and Sobolev spaces were considered), many other functional spaces were used such as Fourier-Herz spaces \cite{cannone2012global,le1997cascades,lei2011global}, Besov spaces \cite{cannone1994ondelettes,lemarie2002recent}, Morrey spaces \cite{giga1989navier,kato1992strong,taylor1992analysis}, the $BMO^{-1}$ space \cite{KochTataru}, etc. For a more complete review of possible functional spaces, see the book \cite{lemarie2018navier}.\\ 

The main feature of this work is to explore some existence results for equations (\ref{NS_Intro}) using as base spaces the Lebesgue spaces of variable exponent $L^{p(\cdot)}$, which -to the best of our knowledge- were not used before in the analysis of the Navier-Stokes equations. These spaces are quite different from the usual Lebesgue spaces $L^p$. Indeed, the parameter $p$ now is a function $p(\cdot):\mathbb{R}^3\longrightarrow [1,+\infty[$. To define the spaces  $L^{p(\cdot)}$ we will proceed as follows: for a measurable function $\vf:\mathbb{R}^3\longrightarrow \mathbb{R}^3$, we consider the \emph{modular function} $\varrho_{p(\cdot)}$ associated to $p(\cdot)$, which is given by the expression
\begin{equation}\label{Def_Modular_Intro}
\varrho_{p(\cdot)}(\vf)=\int_{\mathbb{R}^3}|\vf(x)|^{p(x)}dx.
\end{equation}
We note that if the function $p(\cdot)$ is constant (\textit{i.e.} if $p(\cdot)\equiv p\in [1,+\infty[$)  we obtain the classical Lebesgue spaces and we can derive from the modular function  $\varrho_p$ a norm defined by 
$$\|\vf\|_{L^{p}}=\left(\int_{\mathbb{R}^3}|\vf(x)|^{p}dx\right)^{\frac1p}.$$
However, in the general case where $p(\cdot)$ is a measurable function, it is not possible replace in the previous formula the constant exponent $\frac1p$ outside the integral by $\frac{1}{p(\cdot)}$. To overcome this issue it is classical (see the books \cite{Cruz_Libro} and \cite{Diening_Libro}) to consider the \emph{Luxemburg norm} $\|\cdot\|_{L^{p(\cdot)}}$ associated to the modular function $\varrho_{p(\cdot)}$, which is given by:
\begin{equation}\label{Def_LuxNormLebesgue}
\|\vf\|_{L^{p(\cdot)}}=\inf\{\lambda > 0: \, \varrho_{p(\cdot)}(\vf/\lambda)\leq1\}.
\end{equation}
We then define the Lebesgue spaces of variable exponent $L^{p(\cdot)}(\mathbb{R}^3)$ as the set of all the measurable functions such that the quantity $\|\cdot\|_{L^{p(\cdot)}}$ given above is finite  (for more details we refer the reader to the previously cited books or to the Section \ref{Secc_Notaciones_Presentaciones} below).\\

In our results will thus study mild solutions for the Navier-Stokes system in the framework of the Lebesgue spaces of variable exponent. Indeed, in Theorem \ref{Theoreme_1} below we will consider a particular variant of these functional spaces for the space variable $x\in \mathbb{R}^3$ and we will consider a classical $L^\infty$ space for the time variable $t>0$. In Theorem \ref{Theoreme_2} we will work with a Lebesgue space of variable exponent in the time variable with a usual $L^q$ space in the space variable. In these two different cases we will be able to close the fixed point argument and to obtain mild solutions: these results give a first idea of what can be done in this framework of Lebesgue spaces of variable exponent. Moreover, since the $L^{p(\cdot)}$ spaces can not be easily related to the usual Lebesgue spaces, we expect that these results could offer a different point of view in some applications.\\[1mm]

The complete analysis of the Navier-Stokes equations is a deep and complex field and our aim here is to present a first application of the Lebesgue spaces of variable exponent to this topic. We thus hope that this work will help to attract the attention to these functional spaces.\\

The outline of the article is the following. In Section \ref{Secc_Notaciones_Presentaciones} we first present a small review of the main properties of the spaces $L^{p(\cdot)}$ and then we state our results. Section \ref{Secc_Proof_Existence} is devoted to the proof of the theorems.
%%%%%%%%%%%%%%%%%%%%%%%%%%%%%%%%%%%%%%%%%%%%%%%%%%%
\section{Preliminaries and presentation of results}\label{Secc_Notaciones_Presentaciones}
For $n\geq 1$, let us first consider a function $p:\mathbb{R}^n\longrightarrow [1,+\infty[$, we will say that $p\in \mathcal{P}(\mathbb{R}^n)$ if $p(\cdot)$ is a measurable function and we define $p^-=\underset{x\in \mathbb{R}^n}{\mbox{inf ess}} \; \{p(x)\}$ and $p^+=\underset{x\in \mathbb{R}^n}{\mbox{sup ess}} \; \{p(x)\}$. In order to distinguish between variable and constant exponents, we will always denote exponent functions by $p(\cdot)$, moreover, for the sake of simplicity and to avoid technicalities 
(see \cite[Chapter 3]{Diening_Libro}),
we will always assume here that we have
\begin{equation}\label{BornePminPmax}
1<p^-\leq p^+<+\infty.
\end{equation}
With these exponents we can consider the Luxemburg norm $\|\cdot\|_{L^{p(\cdot)}}$ as defined in  (\ref{Def_LuxNormLebesgue}). These functional spaces $L^{p(\cdot)}(\mathbb{R}^n)$ possess some of the structural properties of normed spaces (they are moreover Banach function spaces) but they also present some very special features. \\

In this setting, the H\"older inequalities have the following version: let $p(\cdot),\,q(\cdot),\,r(\cdot)\in \mathcal{P}(\mathbb{R}^n)$ be functions such that we have the pointwise relationship
$\frac{1}{p(x)}=\frac{1}{q(x)}+\frac{1}{r(x)}$, $x\in \mathbb{R}^n$. Then there exists a constant $C>0$ such that for all $f\in L^{q(\cdot)}(\mathbb{R}^n) $ and $g \in L^{r(\cdot)}(\mathbb{R}^n)$, the pointwise product $fg$ belongs to the space $L^{p(\cdot)}(\mathbb{R}^n)$ and we have the estimate
\begin{equation}\label{Holder_LebesgueVar}
\|f g\|_{L^{p(\cdot)}} \leq C\|f\|_{L^{q(\cdot)}}\|g\|_{L^{r(\cdot)}},
\end{equation}
see \cite[Section 2.4]{Cruz_Libro} or \cite[Section 3.2]{Diening_Libro} for a proof of this fact. This estimate can be easily generalized to vector fields $\vf, \vg :\mathbb{R}^n\longrightarrow \mathbb{R}^n$ and to the product $\vf\cdot \vg$.\\

Note that the quantity $\|\cdot\|_{L^{p(\cdot)}}$ satisfies the \emph{Norm conjugate formula} given in \cite[Corollary 3.2.14]{Diening_Libro}:
\begin{equation}\label{Norm_conjugate_formula}
\|f\|_{L^{p(\cdot)}}\leq \underset{\|g\|_{L^{p'(\cdot)}}\leq 1}{\sup}\int_{\mathbb{R}^n}|f(x)||g(x)|dx\qquad \mbox{with } \frac{1}{p(\cdot)}+\frac{1}{p'(\cdot)}=1.
\end{equation}
%%%%%%%%%%%%%%%%%%%%%%%%%%%%%%%%%%%%%%%%%%%%%%%%%%%
\begin{Remarque}
Note that in the previous notions the space $\mathbb{R}^n$ can be replaced by an interval $[0,T]$. 
\end{Remarque}
%%%%%%%%%%%%%%%%%%%%%%%%%%%%%%%%%%%%%%%%%%%%%%%%%%%
It is crucial to remark now that the convolution product $f\ast g$ is not well adapted to the structure of the $L^{p(\cdot)}$ spaces, in particular the Young inequalities for convolution are not valid anymore (see \cite[Section 3.6]{Diening_Libro}) and thus many of the usual operators that appear in PDEs must be treated very carefully. Note also that Fourier-based methods are not so easy to use as we lack of an alternative for the Plancherel formula.\\

To study the boundedness of such operators we first need to impose some conditions over the functions $p(\cdot)\in \mathcal{P}(\mathbb{R}^n)$: indeed, following \cite[Section 4.1]{Diening_Libro}, we will say that a measurable function $p(\cdot)\in \mathcal{P}(\mathbb{R}^n)$ belongs to the class $\mathcal{P}^{log}(\mathbb{R}^n)$ if we have 
$$\left|\frac{1}{p(x)}-\frac{1}{p(y)}\right|\leq \frac{C}{\log(e+1/|x-y|)}\qquad \mbox{for all } x,y\in \mathbb{R}^n,$$
and if 
$$\left|\frac{1}{p(x)}-\frac{1}{p_\infty}\right| \leq \frac{C}{\log(e+|x|)}\qquad \mbox{for all } x\in \mathbb{R}^n,$$
where 
\begin{equation}
\frac{1}{p_\infty}=\underset{|x|\to +\infty}{\lim}\frac{1}{p(x)}. 
\end{equation}
With the condition $p(\cdot)\in \mathcal{P}^{log}(\mathbb{R}^n)$ we have the following results.
\begin{itemize}
\item For $f:\mathbb{R}^n\longrightarrow \mathbb{R}$ a locally integrable function, the Hardy-Littlewood maximal function $\mathcal{M}$ is given by $\displaystyle{\mathcal{M}(f)(x)=\underset{B \ni x}{\sup } \;\frac{1}{|B|}\int_{B }|f(y)|dy}$ where $B$ is an open ball of $\mathbb{R}^n$. Thus, if $p(\cdot)\in \mathcal{P}^{log}(\mathbb{R}^n)$ we have the estimate
\begin{equation}\label{MaximalFunc_LebesgueVar}
\|\mathcal{M}(f)\|_{L^{p(\cdot)}}\leq C \|f\|_{L^{p(\cdot)}}.
\end{equation}
See \cite[Section 4.3]{Diening_Libro} for a proof of this fact.\\

\item Note also that the usual Riesz transforms $(\mathcal{R}_j)_{1\leq j\leq n}$ defined formally in the Fourier level by $\widehat{\mathcal{R}_j(f)}(\xi)=-\frac{i\xi_j}{|\xi|}\widehat{f}(\xi)$ are also bounded in Lebesgue spaces of variable exponent and we have the inequality
\begin{equation}\label{Riesz_LebesgueVar}
\|\mathcal{R}_j(f)\|_{L^{p(\cdot)}}\leq C \|f\|_{L^{p(\cdot)}},
\end{equation}
with the conditions  $p(\cdot)\in \mathcal{P}^{log}(\mathbb{R}^n)$ and $1<p^-\leq p^+<+\infty$. See 
\cite[Section 6.3]{Diening_Libro}.\\

\item Let us recall now that for $0<\sigma<n$, the Riesz potentials $\mathcal{I}_\sigma$ are defined by 
\begin{equation}\label{Definition_RieszPotential}
\mathcal{I}_\sigma(f)(x)=\int_{\mathbb{R}^n}\frac{|f(y)|}{|x-y|^{n-\sigma}}dy. 
\end{equation}
If $p(\cdot)\in \mathcal{P}^{log}(\mathbb{R}^n)$ and if $0<\sigma<n/p^+$, then, following \cite[Section 6.1]{Diening_Libro}, we have the inequality 
\begin{equation}\label{PotentialRieszVariable0}
\|\mathcal{I}_\sigma(f)\|_{L^{q(\cdot)}}\leq C\|f\|_{L^{p(\cdot)}},\qquad \mbox{with} \quad  \frac{1}{q(\cdot)}=\frac{1}{p(\cdot)}-\frac{\sigma}{n}.
\end{equation}
This estimate introduces a very strong relationship between the parameters $p(\cdot)$ and $q(\cdot)$. In order to obtain some more freedom in the parameters (see Remark \ref{Rem_Riesz_MixedLebesgue} below) we will consider the mixed Lebesgue spaces $\mathcal{L}^{p(\cdot)}_\mathfrak{p}(\mathbb{R}^n)=L^{p(\cdot)}(\mathbb{R}^n)\cap L^{\mathfrak{p}}(\mathbb{R}^n)$ introduced in \cite{chamorro_paper_lpvar}, where $1<\mathfrak{p}<+\infty$ is a constant exponent. These spaces that can be normed by the quantity 
\begin{equation}\label{MixedLebesgue}
\|\cdot\|_{\mathcal{L}^{p(\cdot)}_\mathfrak{p}}=\max\{\|\cdot\|_{L^{p(\cdot)}}, \|\cdot\|_{L^{\mathfrak{p}}}\}.
\end{equation}
With the help of these spaces we have the following result. 
%%%%%%%%%%%%%%%%%%%%%%%%%%%%%%%%%%%%%%%%%%%%%%%%%%%
\begin{Proposition}\label{Proposition_RieszPotential}
Let $1<\mathfrak{p}<+\infty$ be a constant exponent, $p(\cdot)\in \mathcal{P}^{\log}(\mathbb{R}^n)$ a variable exponent and fix a parameter $0<\sigma<\min\{n/p^+, n/\mathfrak{p}\}$. If $f\in \mathcal{L}^{p(\cdot)}_\mathfrak{p}(\mathbb{R}^n)$, then we have the inequality
\begin{equation}\label{Inegalite3}
\|\mathcal{I}_\sigma(f)\|_{L^{\rho(\cdot)}}\leq C \|f\|_{\mathcal{L}^{p(\cdot)}_\mathfrak{p}},
\end{equation}
where the function $\rho(\cdot)$ satisfies the following condition
\begin{equation}\label{Inegalite2Condition}
\rho(\cdot)=\frac{np(\cdot)}{n-s\mathfrak{p}}.
\end{equation}
\end{Proposition}
%%%%%%%%%%%%%%%%%%%%%%%%%%%%%%%%%%%%%%%%%%%%%%%%%%%
See \cite{chamorro_paper_lpvar} for a proof. Note in particular that the index $\mathfrak{p}$ is not to related to $p^-$ or $p^+$ nor to $p(\cdot)$ and this inequality gives more flexibility in the indexes than the conditions of the estimate (\ref{PotentialRieszVariable0}).
\begin{Remarque}\label{Rem_Holder_Mixed_Lebesgue_Var}
Note that, by construction, the mixed spaces $\mathcal{L}^{p(\cdot)}_\mathfrak{p}$ inherit the properties of the spaces $L^{p(\cdot)}$ and $L^{\mathfrak{p}}$. In particular we have the H\"older inequality $\|\varphi\|_{\mathcal{L}^{p(\cdot)}_\mathfrak{p}}\leq \|\varphi\|_{\mathcal{L}^{q(\cdot)}_\mathfrak{q}}\|\varphi\|_{\mathcal{L}^{r(\cdot)}_\mathfrak{r}}$ with $\frac{1}{p(\cdot)}=\frac{1}{q(\cdot)}+\frac{1}{r(\cdot)}$ and $\frac{1}{\mathfrak{p}}=\frac{1}{\mathfrak{q}}+\frac{1}{\mathfrak{r}}$ and of course the Riesz transforms are also bounded in these spaces.\\
\end{Remarque}
For more details on the Lebesgue spaces of variable exponent, on their inner structure as well as many other properties, see the books \cite{Cruz_Libro} and \cite{Diening_Libro}.\\
\end{itemize}
With these preliminaries, which were presented for the sake of generality in $\mathbb{R}^n$, we can state our first result about the existence of mild solutions for the 3D Navier-Stokes system:
%%%%%%%%%%%%%%%%%%%%%%%%%%%%%%%%%%%%%%%%%%%%%%%%%%% 
\begin{Theoreme}[Global Mild Solutions]\label{Theoreme_1}
Consider a variable exponent $p(\cdot)\in \mathcal{P}^{\log}(\mathbb{R}^3)$, a divergence free initial data $\vu_0\in \mathcal{L}^{p(\cdot)}_{3}(\mathbb{R}^3)$ (as defined in (\ref{MixedLebesgue}) above) and let $\vf$ be a divergence free external force such that $\vf=div(\mathbb{F})$ where $\mathbb{F}$ is a tensor such that $ \mathcal{L}^{\frac{p(\cdot)}{2}}_{\frac{3}{2}}(\mathbb{R}^3, L^\infty([0,T[))$. If the quantity $\|\vu_0\|_{\mathcal{L}^{p(\cdot)}_{3}}+\|\mathbb{F}\|_{\mathcal{L}^{\frac{p(\cdot)}{2}}_{\frac{3}{2},x}(L^\infty_t)}$ is small, then the Navier-Stokes equations (\ref{NS_Intro}) admits a unique, global mild solution in the space $\mathcal{L}^{p(\cdot)}_{3}\left(\mathbb{R}^3,L^\infty([0,T[)\right)$.
\end{Theoreme}
%%%%%%%%%%%%%%%%%%%%%%%%%%%%%%%%%%%%%%%%%%%%%%%%%%%
Some remarks are in order here. First note that as long as the initial data and the external force are small enough, then we can obtain a unique, global solution for the Navier-Stokes equations (\ref{NS_Intro}) in the framework of Lebesgue spaces of variable exponent. This seems to be, to the best of our knowledge, the first application of this type of spaces in the analysis of these equations. Let us mention for the sake of completness that steady versions of some PDEs from fluid mechanics were studied in the Section 14.4 of the book \cite{Diening_Libro} but evolution problems, as the one considered here, require a different treatment. Remark also that global mild solution are frequently associated to small data, however the study of large initial data that could generate global mild solution is a completely (and hard) open problem. Finally, let us point out that the use of the mixed Lebesgue spaces of variable exponent $\mathcal{L}^{p(\cdot)}_{\mathfrak{p}}$ is essentially technical and it is driven by the lack of flexibility of the indexes that intervene in the boundedness of the Riesz transforms. See the Remark \ref{Rem_Riesz_MixedLebesgue} below for more details in this particular point.\\

In the previous result, we have first analyzed the behavior of the solutions in the time variable (in a $L^\infty$ norm) and then we studied the information in the space variable in a mixed Lebesgue space of variable exponent. In our second theorem we will proceed in a different fashion: 
%%%%%%%%%%%%%%%%%%%%%%%%%%%%%%%%%%%%%%%%%%%%%%%%%%% 
\begin{Theoreme}[Local Mild Solutions ]\label{Theoreme_2}
Let $p(\cdot )\in \mPl(\Rt)$ with $p^->2$ and fix an index $q>3$ by the relationship $\frac{2}{p(\cdot)} + \frac{3}{q} < 1$. If $\vf\in L^{p(\cdot)} \left( [0,+\infty[,  L^{q} (\Rt) \right)$ is an exterior force and if $\vu_0\in L^{q} (\Rt)$ is an initial data such that $div(\vu_0)=0$, then there exist a time $0<T<+\infty$ and an unique local in time mild solution of the Navier-Stokes equations (\ref{NS_Intro}) in the space $L^{p(\cdot)} \left( [0,T], L^{q} (\Rt) \right)$.
\end{Theoreme}
%%%%%%%%%%%%%%%%%%%%%%%%%%%%%%%%%%%%%%%%%%%%%%%%%%% 
Note that, if we compare this result with Theorem \ref{Theoreme_1}, we changed here the order of the variables: we first measure the information in the space variable and then we consider the information in the time variable. Although it is a slightly more popular way to construct mild solutions for evolutive PDEs, we face here the inexorable problem of the time of existence of such solutions: large initial conditions can be considered but then the time of existence will be very small. \\

Note moreover that, if the exponent $p$ is constant, the condition $\frac{2}{p} + \frac{3}{q} < 1$ is common in the analysis of the Navier-Stokes system, however we should expect that the additional freedom given by the variable exponent $p(\cdot)$ would be used to deepen the study of these equations.

%%%%%%%%%%%%%%%%%%%%%%%%%%%%%%%%%%%%%%%%%%%%%%%%%%%
\section{Mild solutions in Lebesgue Spaces of variable exponent}\label{Secc_Proof_Existence}
We present here a first general approach to mild solutions for the Navier-Stokes equations (\ref{NS_Intro}) in the setting of Lebesgue spaces of variable exponent. These mild solutions are obtained via the following classical result:
%%%%%%%%%%%%%%%%%%%%%%%%%%%%%%%%%%%%%%%%%%%%%%%%%%%
\begin{Theoreme}[Banach-Picard contraction principle]\label{BP_principle}
Let $(E,\|\cdot \|_E )$ a Banach space and consider $B: E \times E \longrightarrow E$ a bounded bilinear application: 
$$\|B(e,e)\|_{E}\leq C_B\|e\|_E\|e\|_E.$$
Given $e_0\in E$  such that $\|e_0\|_E \leq \delta$ with $0<\delta < \frac{1}{4C_B}$, then the equation 
$$e = e_0 -  B(e,e),$$
admits an unique solution $e\in E$ which satsifies $\| e \|_E \leq 2 \delta$.
\end{Theoreme}
%%%%%%%%%%%%%%%%%%%%%%%%%%%%%%%%%%%%%%%%%%%%%%%%%%%
In order to apply this result to the Navier-Stokes equations (\ref{NS_Intro}) we need to get rid of the pressure $P$ and for this we apply to this system the Leray projector $\mathbb{P}$ defined by\footnote{Recall that the Leray projector $\mathbb{P}$ can also be defined in terms of the Riesz transforms: $\mathbb{P}(\vec{\varphi})=(Id_{3\times 3}-\vec{R}\otimes\vec{R})(\vec{\varphi})$ where $\vec{R}=(R_1, R_2, R_3)$ and $R_j$ is the $j$-th Riesz transform. Thus, as long as the Riesz transforms are bounded in a functional space $E$, then the Leray projector $\mathbb{P}$ is bounded in $E$ and we have $\|\mathbb{P}(\vec{\varphi})\|_E\leq C\|\vec{\varphi}\|_E$.} $\mathbb{P}(\vec{\varphi})=\vec{\varphi}+\vn \frac{1}{(-\Delta)}(\vn\cdot \vec{\varphi})$. Recall that we thus have $\mathbb{P}(\vn P)\equiv 0$ and if a vector field is divergence free we have the identity $\mathbb{P}(\vec{\varphi})= \vec{\varphi}$. Since $\vu$ and $\vf$ are divergence free, we obtain the equation
$$
\begin{cases}
\partial_t\vu=\Delta \vu-\mathbb{P}(div(\vu \otimes \vu))+\vf, \quad div(\vu)=0,\\[3mm]
\vu(0,x)=\vu_0(x), \qquad x\in \mathbb{R}^3.
\end{cases}
$$
Now, due to the Dumahel formula, we can write this equation in the following form 
\begin{equation}\label{NS_Integral}
\vu(t,x)=\mathfrak{g}_t\ast \vu_0(x)+\int_{0}^t\mathfrak{g}_{t-s}\ast\vf(s, x)ds-\int_{0}^t\mathfrak{g}_{t-s}\ast \mathbb{P}(div(\vu \otimes \vu))(s, x)ds,
\end{equation}
where $\mathfrak{g}_t$ is the usual gaussian heat kernel.  \\

%%%%%%%%%%%%%%%%%%%%%%%%%%%%%%%%%%%%%%%%%%%%%%%%%%%
\subsection{Proof of Theorem \ref{Theoreme_1}}
The integral equation (\ref{NS_Integral}) above is now of the form considered in Theorem \ref{BP_principle}: indeed, it is enough to set the quantity 
$$\mathfrak{g}_t\ast \vu_0(x)+\int_{0}^t\mathfrak{g}_{t-s}\ast\vf(s, x)ds,$$
as the term $e_0$ and to analyse the bilinear application 
$$B(\vu, \vu)=\int_{0}^t\mathfrak{g}_{t-s}\ast \mathbb{P}(div(\vu \otimes \vu))(s, x)ds.$$
We will thus look for a mild solution of equation \eqref{NS_Integral} in the functional space 
$\mathcal{E}=\mathcal{L}^{p(\cdot)}_{3}(\mathbb{R}^3,L^\infty([0,T[))$ endowed with the norm
\begin{equation}\label{Norm_PointFixe}
\|\cdot\|_{\mathcal{E}}=\max\{\|\cdot\|_{L^{p(\cdot)}_x(L^\infty_t)}, \|\cdot\|_{L^{3}_x(L^\infty_t)}\},
\end{equation}
where $p(\cdot)\in \mathcal{P}^{\log}(\mathbb{R}^3)$ and $1<p^-\leq p^+<+\infty$. Now, following the argument presented in Theorem \ref{BP_principle}, we will prove the following estimates
\begin{equation}\label{Control_Uo}
\|\mathfrak{g}_t\ast \vu_0\|_{\mathcal{E}}\leq C\|\vu_0\|_{\mathcal{L}^{p(\cdot)}_{3}},
\end{equation}
\begin{equation}\label{Control_f_preambulo}
\left\|\int_{0}^t\mathfrak{g}_{t-s}\ast \vf(\cdot, \cdot)ds\right\|_{\mathcal{E}}
\leq C\|\mathbb{F}\|_{\mathcal{L}^{\frac{p(\cdot)}{2}}_{\frac{3}{2},x}(L^\infty_t)}.
\end{equation}
and 
\begin{equation}\label{Control_Nolineal}
\left\|\int_{0}^t\mathfrak{g}_{t-s}\ast \mathbb{P}(div(\vu \otimes \vu))(\cdot, \cdot)ds\right\|_{\mathcal{E}}\leq C_B\|\vu\|_{\mathcal{E}}\|\vv\|_{\mathcal{E}}.
\end{equation}
Thus, if we have the condition
$$C\big(\|\vu_0\|_{\mathcal{L}^{p(\cdot)}_{3}}+\|\mathbb{F}\|_{\mathcal{L}^{\frac{p(\cdot)}{2}}_{\frac{3}{2},x}(L^\infty_t)}\big)<\frac{1}{4C_B},$$
then we can obtain a unique mild solution for the system (\ref{NS_Integral}).\\

Now, we will deduce each one of the previous estimates (\ref{Control_Uo}),  (\ref{Control_f_preambulo}) and (\ref{Control_Nolineal}). 
\begin{itemize}
\item First we study the quantity (\ref{Control_Uo}) and for this we recall a classical lemma (see \cite[Lemma 7.4, Section 7.7]{lemarie2018navier}):
%%%%%%%%%%%%%%%%%%%%%%%%%%%%%%%%%%%%%%%%%%%%%%%%%%%
\begin{Lemme}\label{lemme_conv_maximal}
If $\varphi$ is a radially decreasing function on $\mathbb{R}^3$ and $\vf$ is a locally integrable function, then 
\begin{equation*}
|(\varphi\ast\vf)(x)| \leq \|\varphi\|_{L^1}  \mathcal{M} (\vf)(x),
\end{equation*}
where $\mathcal{M}$ is the Hardy-Littlewood maximal function.
\end{Lemme}
%%%%%%%%%%%%%%%%%%%%%%%%%%%%%%%%%%%%%%%%%%%%%%%%%%%
Since the heat kernel $\mathfrak{g}_t$ is a radially decreasing function and since $\vu_0$ is a locally integrable function, using the previous lemma we obtain the control
$$\|\mathfrak{g}_t\ast \vu_0(x)\|_{L^\infty_t}\leq C\mathcal{M}(\vu_0)(x),$$
and we can write
\begin{eqnarray*}
\|\mathfrak{g}_t\ast \vu_0\|_{\mathcal{E}}&=&\max\{\|\mathfrak{g}_t\ast \vu_0\|_{L^{p(\cdot)}_x(L^\infty_t)}, \|\mathfrak{g}_t\ast \vu_0\|_{L^{3}_x(L^\infty_t)}\}\\
&\leq & C\max\{\|\mathcal{M}(\vu_0)\|_{L^{p(\cdot)}}, \|\mathcal{M}(\vu_0)\|_{L^{3}}\}.
\end{eqnarray*}
Since $p(\cdot)\in \mathcal{P}^{log}(\mathbb{R}^3)$,  the estimate  (\ref{MaximalFunc_LebesgueVar}) implies that the maximal function $\mathcal{M}$ is bounded in the Lebesgue space $L^{p(\cdot)}(\mathbb{R}^3)$ (this operator is also bounded in $L^3$). Then we obtain
$$\|\mathfrak{g}_t\ast \vu_0\|_{\mathcal{E}}\leq C\max\big\{\|\vu_0\|_{L^{p(\cdot)}},\|\vu_0\|_{L^{3}}\big\}\leq  C\|\vu_0\|_{\mathcal{L}^{p(\cdot)}_3},$$
which is the announced control (\ref{Control_Uo}).

\item To study the inequality \eqref{Control_f_preambulo} we proceed as follows: since $\vf=div(\mathbb{F})$ we can write
$$\left|\int_{0}^t\mathfrak{g}_{t-s}\ast \vf(s,x)ds\right|\leq C\int_{0}^t\int_{\mathbb{R}^3}|\vn\mathfrak{g}_{t-s}(x-y)| |\mathbb{F}(s,y)|dyds,$$
and due to the decay properties of the heat kernel we obtain
\begin{eqnarray*}
\left|\int_{0}^t\mathfrak{g}_{t-s}\ast \vf(s,x)ds\right|&\leq &C\int_{0}^t\int_{\mathbb{R}^3}\frac{1}{|t-s|^2+|x-y|^4} |\mathbb{F}(s,y)|dyds\\
&\leq &C\int_{\mathbb{R}^3}\int_{0}^t\frac{1}{|t-s|^2+|x-y|^4} |\mathbb{F}(s,y)|dsdy,
\end{eqnarray*}
where we have applied the Fubini Theorem. Moreover, since $\mathbb{F}\in \mathcal{L}^{\frac{p(\cdot)}{2}}_{\frac{3}{2},x}(L^\infty_t)$, one has
$$\left|\int_{0}^t\mathfrak{g}_{t-s}\ast \vf(s,x)ds\right|\leq C\int_{\mathbb{R}^3}\int_{0}^t\frac{1}{|t-s|^2+|x-y|^4} ds \|\mathbb{F}(\cdot,y)\|_{L^\infty_t}dy,$$
and after an integration with respect to the time variable, it comes
$$\left|\int_{0}^t\mathfrak{g}_{t-s}\ast \vf(s,x)ds\right|\leq C\int_{\mathbb{R}^3}\frac{1}{|x-y|^2} \|\mathbb{F}(\cdot,y)\|_{L^\infty_t}dy=\mathcal{I}_1(\|\mathbb{F}(\cdot,\cdot)\|_{L^\infty_t})(y),$$
where $\mathcal{I}_1$ is the Riesz potential defined in (\ref{Definition_RieszPotential}). We thus obtain
$$\left\|\int_{0}^t\mathfrak{g}_{t-s}\ast \vf(s,x)ds\right\|_{L^\infty_t}\leq C\mathcal{I}_1(\|\mathbb{F}(\cdot,\cdot)\|_{L^\infty_t})(y).$$
Now, in order to reconstruct the $\mathcal{L}^{p(\cdot)}_3$ norm given in (\ref{Norm_PointFixe}), from the previous estimate we write
\begin{eqnarray*}
\left\|\int_{0}^t\mathfrak{g}_{t-s}\ast \vf(s,x)ds\right\|_{L^{p(\cdot)}_x(L^\infty_t)}&\leq &C\left\|\mathcal{I}_1(\|\mathbb{F}(\cdot,\cdot)\|_{L^\infty_t})(\cdot)\right\|_{L^{p(\cdot)}_x}\\
\left\|\int_{0}^t\mathfrak{g}_{t-s}\ast \vf(s,x)ds\right\|_{L^{3}_x(L^\infty_t)}&\leq &C\left\|\mathcal{I}_1(\|\mathbb{F}(\cdot,\cdot)\|_{L^\infty_t})(\cdot)\right\|_{L^{3}_x}.
\end{eqnarray*}
Thus by the Proposition \ref{Proposition_RieszPotential} with $\sigma=1$, $\mathfrak{p}=\frac{3}{2}$ and $\rho(\cdot)=p(\cdot)$, we have the following estimate:
$$\left\|\mathcal{I}_1(\|\mathbb{F}(\cdot,\cdot)\|_{L^\infty_t})(\cdot)\right\|_{L^{p(\cdot)}_x}\leq C\left\|\|\mathbb{F}(\cdot,\cdot)\|_{L^\infty_t}\right\|_{\mathcal{L}^{\frac{p(\cdot)}{2}}_{\frac{3}{2},x}}=\|\mathbb{F}\|_{\mathcal{L}^{\frac{p(\cdot)}{2}}_{\frac{3}{2},x}(L^\infty_t)}.$$
Moreover, by the boundedness properties of the Riesz potentials in the usual Lebesgue spaces $L^3$ we obtain 
$$\left\|\mathcal{I}_1(\|\mathbb{F}(\cdot,\cdot)\|_{L^\infty_t})(\cdot)\right\|_{L^{3}_x}\leq C\left\|\|\mathbb{F}\|_{L^\infty_t}\right\|_{L^{\frac{3}{2}}_x}=\|\mathbb{F}\|_{L^{\frac{3}{2}}_x(L^\infty_t)}.$$
With these two estimates at hand, and following the definition of the norm $\|\cdot\|_{\mathcal{E}_T}$ given in (\ref{Norm_PointFixe}), we finally obtain
$$\left\|\int_{0}^t\mathfrak{g}_{t-s}\ast \vf(s,x)ds\right\|_{\mathcal{E}}=\left\|\int_{0}^t\mathfrak{g}_{t-s}\ast \vf(s,x)ds\right\|_{\mathcal{L}^{p(\cdot)}_{3,x}(L^\infty_t)}\leq C\|\mathbb{F}\|_{\mathcal{L}^{\frac{p(\cdot)}{2}}_{\frac{3}{2},x}(L^\infty_t)}<+\infty.$$

\item In order to establish the estimate (\ref{Control_Nolineal}), we first remark that, due to the properties of the Leray projector, we have the identities
$$\int_{0}^t\mathfrak{g}_{t-s}\ast \mathbb{P}(div(\vu \otimes \vv))ds=\int_{0}^t\mathbb{P}\big(\mathfrak{g}_{t-s}\ast div(\vu \otimes \vu)\big)ds=\mathbb{P}\left(\int_{0}^t\mathfrak{g}_{t-s}\ast div(\vu \otimes \vu)ds\right),$$
and as before, by the decay properties of the heat kernel, we obtain
\begin{eqnarray*}
\left|\int_{0}^t\mathfrak{g}_{t-s}\ast div(\vu \otimes \vu)(s,x)ds\right|&\leq &C\int_{0}^t\int_{\mathbb{R}^3}|\vn \mathfrak{g}_{t-s}(x-y)|\vu(s,y)| |\vu(s,y)|dyds\\
&\leq & C\int_{\mathbb{R}^3}\int_{0}^t\frac{1}{|t-s|^2+|x-y|^4} |\vu(s,y)| |\vu(s,y)|dyds.
\end{eqnarray*}
Note that we have $|\vu(t,x)|\leq \|\vu(\cdot,x)\|_{L^\infty_t}$ and we can thus write
$$\left|\int_{0}^t\mathfrak{g}_{t-s}\ast div(\vu \otimes \vu)(s,x)ds\right|\leq C\int_{\mathbb{R}^3}\left(\int_{0}^t \frac{1}{|t-s|^2+|x-y|^4}ds\right)\|\vu(\cdot,x)\|_{L^\infty_t}\|\vu(\cdot,x)\|_{L^\infty_t}dy,$$
from which we deduce the estimate
$$\left|\int_{0}^t\mathfrak{g}_{t-s}\ast div(\vu \otimes \vu)(s,x)ds\right|\leq C\int_{\mathbb{R}^3}\frac{1}{|x-y|^2}\|\vu(\cdot,x)\|_{L^\infty_t}\|\vu(\cdot,x)\|_{L^\infty_t}dy.$$
Remark now that, using the definition of the Riesz potentials given in 
(\ref{Definition_RieszPotential}), the last term above can be written in the following form:
$$\left|\int_{0}^t\mathfrak{g}_{t-s}\ast div(\vu \otimes \vu)(s,x)ds\right|\leq C\mathcal{I}_1\big(\|\vu\|_{L^\infty_t}\|\vu\|_{L^\infty_t}\big)(x),$$
and we obtain the estimate
$$\mathbb{P}\left(\int_{0}^t\mathfrak{g}_{t-s}\ast div(\vu \otimes \vu)(s,x)ds\right)\leq C\mathbb{P}\left(\mathcal{I}_1\big(\|\vu\|_{L^\infty_t}\|\vu\|_{L^\infty_t}\big)(x)\right).$$
In order to reconstruct the norm $\|\cdot\|_{\mathcal{E}_T}$ given in (\ref{Norm_PointFixe}) we write:
\begin{eqnarray*}
\left\|\mathbb{P}\left(\int_{0}^t\mathfrak{g}_{t-s}\ast div(\vu \otimes \vu)ds\right)\right\|_{L^{p(\cdot)}_x(L^\infty_t)}&\leq &C\|\mathcal{I}_1(\|\vu\|_{L^\infty_t}\|\vu\|_{L^\infty_t})\|_{L^{p(\cdot)}_x(L^\infty_t)}\\
\left\|\mathbb{P}\left(\int_{0}^t\mathfrak{g}_{t-s}\ast div(\vu \otimes \vu)ds\right)\right\|_{L^{3}_x(L^\infty_t)}&\leq &C\|\mathcal{I}_1(\|\vu\|_{L^\infty_t}\|\vu\|_{L^\infty_t})\|_{L^{3}_x(L^\infty_t)},
\end{eqnarray*}
where we have used the fact that the Leray projector $\mathbb{P}$ is bounded in the Lebesgue space $L^3$ as well as in the Lebesgue space of variable exponent $L^{p(\cdot)}$ (since the Riesz transforms are bounded in such spaces). At this point we recall that the Riesz potential $\mathcal{I}_1$ satisfies $\|\mathcal{I}_1(\varphi)\|_{L^3}\leq C\|\varphi\|_{L^{\frac{3}{2}}}$ and following the Proposition \ref{Proposition_RieszPotential} we have $\|\mathcal{I}_1(\varphi)\|_{L^{p(\cdot)}}\leq C \|\varphi\|_{\mathcal{L}^{\frac{p(\cdot)}{2}}_{\frac{3}{2}}}$, so we can write, using the H\"older inequalities (see the Remark \ref{Rem_Holder_Mixed_Lebesgue_Var} above):
\begin{eqnarray*}
\left\|\mathbb{P}\left(\int_{0}^t\mathfrak{g}_{t-s}\ast div(\vu \otimes \vu)ds\right)\right\|_{L^{p(\cdot)}_x(L^\infty_t)}&\leq & C \left\|\|\vu\|_{L^\infty_t}\|\vu\|_{L^\infty_t}\right\|_{\mathcal{L}^{\frac{p(\cdot)}{2}}_{\frac{3}{2}}}\leq  C\|\vu\|_{\mathcal{L}^{p(\cdot)}_{3,x}(L^\infty_x)}\|\vu\|_{\mathcal{L}^{p(\cdot)}_{3,x}(L^\infty_x)} \\
\left\|\mathbb{P}\left(\int_{0}^t\mathfrak{g}_{t-s}\ast div(\vu \otimes \vu)ds\right)\right\|_{L^{3}_x(L^\infty_t)}&\leq &C\left\|\|\vu\|_{L^\infty_t}\|\vu\|_{L^\infty_t}\right\|_{L^{\frac{3}{2}}_x}\leq C\|\vu\|_{L^{3}_x(L^\infty_t)}\|\vu\|_{L^{3}_x(L^\infty_t)}.
\end{eqnarray*}
\begin{Remarque}\label{Rem_Riesz_MixedLebesgue}
It is worth noting here that it is the use of the mixed Lebesgue spaces of variable exponent $\mathcal{L}^{p(\cdot)}_{3}$ that allows us to obtain the inequality $\|\mathcal{I}_1(\varphi)\|_{L^{p(\cdot)}}\leq \|\varphi\|_{\mathcal{L}^{\frac{p(\cdot)}{2}}_{\frac{3}{2}}}$ from which we will deduce the wished control by applying the H\"older inequalities. Indeed, if we only use the spaces $L^{p(\cdot)}$ and we look for an estimate of the form $\|\mathcal{I}_1(\varphi)\|_{L^{p(\cdot)}}\leq \|\varphi\|_{L^{\frac{p(\cdot)}{2}}}$, then the relationship (\ref{PotentialRieszVariable0}) will  force us to have $p(\cdot)\equiv 3$ and this will cancel the interest of using  Lebesgue spaces of variable exponent.
\end{Remarque}
With all these estimates, and using the definition of the norm $\|\cdot\|_{\mathcal{E}}$ given in (\ref{Norm_PointFixe}) we finally can write
$$\left\|\mathbb{P}\left(\int_{0}^t\mathfrak{g}_{t-s}\ast div(\vu \otimes \vu)ds\right)\right\|_{\mathcal{E}}\leq C_B\|\vu\|_{\mathcal{E}}\|\vu\|_{\mathcal{E}}.$$
\end{itemize}
We have established the estimates (\ref{Control_Uo}), (\ref{Control_f_preambulo}) and (\ref{Control_Nolineal}), thus in order to apply the fixed-point Theorem \ref{BP_principle}  we need the condition 
$$C\big(\|\vu_0\|_{\mathcal{L}^{p(\cdot)}_{3}}+\|\mathbb{F}\|_{\mathcal{L}^{\frac{p(\cdot)}{2}}_{\frac{3}{2},x}(L^\infty_t)}\big)<\frac{1}{4C_B},$$
which is fulfilled by hypotheses as the initial data and the external force are assumed to be small, thus the Theorem \ref{Theoreme_1} is proven. \hfill $\blacksquare$

%%%%%%%%%%%%%%%%%%%%%%%%%%%%%%%%%%%%%%%%%%%%%%%%%%%
\subsection{Proof of Theorem \ref{Theoreme_2}} 
As we are interested in mild solutions for the integral problem (\ref{NS_Integral}), we will follow the general steps given in the Theorem \ref{BP_principle} and here we will consider the following functional space 
$$E_{T}=L^{p(\cdot)} \left( [0,T], L^{q} (\Rt) \right),$$ 
with $p(\cdot)\in \mathcal{P}^{log}([0,+\infty[)$, for some $0<T<+\infty$ to be determined later. This space is endowed with the following Luxemburg norm
\begin{equation}\label{Norm_PointFixe_2}
\| \vec{\varphi}\|_{E_T}=\inf\left\{\lambda > 0: \,\int_0^{T}\left|\frac{ \|
\vec{\varphi} (t,\cdot)\|_{L^{q}}}{\lambda}\right|^{p(t)} dt \leq1\right\}.
\end{equation}
\quad\\[5mm]
Thus, in order to apply Theorem \ref{BP_principle}, just as before we first need to establish the following estimates:
\begin{equation}\label{Control_Uo_Lplq}
\|\mathfrak{g}_t * \vu_0  \|_{ E_T}\leq C_1\| \vu_0  \|_{L^{q}},
\end{equation} 
\begin{equation} \label{Control_force_LpLq}
\left\|\int_0^t\mathfrak{g}_{t-s} *\vf(s, \cdot) ds\right\|_{E_T} \leq C_2\|\vf \|_{L^1_t(L^q_x)},
\end{equation} 
and
\begin{equation}\label{Control_Nolineal_LpLq}
\left\|\int_{0}^t\mathfrak{g}_{t-s}\ast \mathbb{P}(div(\vu \otimes \vu))(s, \cdot)ds\right\|_{E_T}  \leq 
C_B\|\vu\|_{E_T}\|\vu\|_{E_T}.
\end{equation}
Note that the constants $C_1, C_2, C_B$ will depend on the time variable as we shall see later on. Each one of these estimates will be studied separately. 
\begin{itemize}
\item To deduce the inequality \eqref{Control_Uo_Lplq}, we simply write
$\|\mathfrak{g}_t * \vu_0  \|_{L^{q} (\Rt)} \leq \|\mathfrak{g}_t   \|_{L^{1} (\Rt)}\| \vu_0  \|_{L^{q} (\Rt)} = \| \vu_0  \|_{L^{q} (\Rt)}$, and to continue, we will need the following result borrowed from the book \cite{Diening_Libro} (Lemma 3.2.12, Section 3.2).
%%%%%%%%%%%%%%%%%%%%%%%%%%%%%%%%%%%%%%%%%%%%%%%%%%%
\begin{Lemme}\label{Lemma_subset}
Let $p(\cdot)\in \mathcal{P}([0,+\infty[)$ with $1<p^-\leq p^+<+\infty$. Then 
$$\frac{1}{C} \min\left\{T^{\frac{1}{p^{-}}}, T^{\frac{1}{p^{+}}}\right\}\leq\|1\|_{L^{p(\cdot)}([0,T])} \leq
C \max\left\{T^{\frac{1}{p^{-}}}, T^{\frac{1}{p^{+}}}\right\}.$$
\end{Lemme}
%%%%%%%%%%%%%%%%%%%%%%%%%%%%%%%%%%%%%%%%%%%%%%%%%%%
Thus, with this result at hand, taking the $L^{p(\cdot)}$ norm in the time variable we have
\begin{equation}\label{Estimation_DonneeInitiale}
\|\mathfrak{g}_t * \vu_0  \|_{L^{p(\cdot)} _tL^{q}_x }\leq C\|\vu_0\|_{L^{q} (\Rt)}\|1\|_{L^{p(\cdot)}([0,T])} \leq C \|\vu_0\|_{L^{q} (\Rt)}\max\left\{T^{\frac{1}{p^{-}}}, T^{\frac{1}{p^{+}}}\right\}.
\end{equation}

\item Let us now prove the estimate \eqref{Control_force_LpLq}. We start by the usual $L^q$-norm in the space variable and we obtain
$$\left\|\int_0^t\mathfrak{g}_{t-s} \ast \vf (s, \cdot)ds\right\|_{L^{q}} \leq  \int_0^t\|\mathfrak{g}_{t-s}\|_{L^1} \| \vf(s, \cdot) \|_{L^{q}}ds \leq C\| \vf\|_{L^1_t(L^{q}_x)}.$$
Proceeding as in the previous lines, \emph{i.e.} taking the $L^{p(\cdot)}$-norm in the time variable and using the Lemma \ref{Lemma_subset} above we can write
\begin{eqnarray}
\left\|\int_0^t\mathfrak{g}_{t-s} \ast \vf(s, \cdot) ds\right\|_{L^{p(\cdot)}_t(L^{q}_x)}&\leq &C\left\|\| \vf  \|_{L^1_t(L^{q}_x)}\right\|_{L^{p(\cdot)}_t}\leq C\| \vf\|_{L^1_t(L^{q}_x)}\|1\|_{L^{p(\cdot)}([0,T])}\notag\\
&\leq &C\| \vf\|_{L^1_t(L^{q}_x)}\max\left\{T^{\frac{1}{p^{-}}}, T^{\frac{1}{p^{+}}}\right\}.\label{Estimation_ForceExterieure}
\end{eqnarray}

\item Let us study now the inequality \eqref{Control_Nolineal_LpLq}. We first take the $L^q$-norm in the space variable to obtain
$$\left\|\int_{0}^t\mathfrak{g}_{t-s}\ast \mathbb{P}(div(\vu \otimes \vu))(s, \cdot)ds
\right\|_{L^q}\leq C\int_{0}^t\left\|\mathfrak{g}_{t-s}\ast \mathbb{P}(div(\vu \otimes \vu))(s, \cdot)
\right\|_{L^q}ds,$$
by the properties of the Leray projector $\mathbb{P}$ and since it is a bounded operator in Lebesgue spaces $L^q$ for $1<q<+\infty$, we can write
$$\left\|\int_{0}^t\mathfrak{g}_{t-s}\ast \mathbb{P}(div(\vu \otimes \vu))(s, \cdot)ds
\right\|_{L^q}\leq C\int_0^t \left\|\vn \mathfrak{g}_{t-s} * \vu\otimes\vu(s, \cdot)\right\|_{L^q} ds.$$
By the Young inequalities with $1+\frac{1}{q}=\frac{q-1}{q}+\frac{2}{q}$ we obtain
\begin{eqnarray*}
\int_0^t \left\|\vn \mathfrak{g}_{t-s} * \vu\otimes\vu(s, \cdot)\right\|_{L^q} ds&\leq &\int_0^t \|\vn \mathfrak{g}_{t-s} \|_{L^{\frac{q}{q-1}}}\left\|\vu\otimes\vu(s, \cdot)\right\|_{L^{\frac{q}{2}}} ds\\
&\leq&C\int_0^t \frac{1}{(t-s)^{\frac{1}{2}+\frac{3}{2q}}}\|\vu(s,\cdot)\|_{L^{q}}\|\vu(s,\cdot)\|_{L^{q}} ds.
\end{eqnarray*}
We will take now the $L^{p(\cdot)}$-norm in the time variable and for this, we will proceed by duality. Thus, for $p'(\cdot)$ the associated conjugate exponent (\emph{i.e.} $\frac{1}{p(\cdot)}+\frac{1}{p'(\cdot)}=1$) we have, by the norm conjugate formula (\ref{Norm_conjugate_formula}), the inequalities:
\begin{eqnarray}
\left\|\int_{0}^t\mathfrak{g}_{t-s}\ast \mathbb{P}(div(\vu \otimes \vu))(s, \cdot)ds
\right\|_{L_t^{p(\cdot)}(L_x^q)}&\leq&\left\|\int_0^t \frac{1}{(t-s)^{\frac{1}{2}+\frac{3}{2q}}}\left\|\vu(s,\cdot)\right\|_{L^{q}}^2ds\right\|_{L_t^{p(\cdot)}([0,T])}\label{Identite_NormeDualite}\\
&\leq &\sup_{\| \psi \|_{L^{p'(\cdot)}}\leq 1} \int_0^{T} \int_0^t \frac{|\psi(t)|}{|t-s|^{\frac{1}{2}+\frac{3}{2q}}}\left\|\vu(s,\cdot)\right\|_{L^{q}}^2ds \,dt.\notag
\end{eqnarray}
We write now, by the Fubini Theorem:
$$\sup_{\| \psi \|_{L^{p'(\cdot)}}\leq 1}\int_0^{T} \int_0^t \frac{|\psi(t)|}{|t-s|^{\frac{1}{2}+\frac{3}{2q}}}\left\|\vu(s,\cdot)\right\|_{L^{q}}^2ds \,dt
=\sup_{\| \psi \|_{L^{p'(\cdot)}}\leq 1}\int_0^T\int_0^T\frac{1_{ \{ 0<s<t \} } |\psi(t)| }{|t-s|^{\frac{1}{2} + \frac{3}{2q} }}dt\|\vu(s,\cdot)
\|^2_{L^q}ds,$$
and extending the function $\psi(t)$ by zero if $t<0$ and if $t>T$, we can write
\begin{eqnarray*}
\sup_{\| \psi \|_{L^{p'(\cdot)}}\leq 1}\int_0^{T} \int_0^t \frac{|\psi(t)|}{|t-s|^{\frac{1}{2}+\frac{3}{2q}}}\left\|\vu(s,\cdot)\right\|_{L^{q}}^2ds \,dt
=\sup_{\| \psi \|_{L^{p'(\cdot)}}\leq 1}\int_0^T\left(\int_{-\infty}^{+\infty}\frac{|\psi(t)| }{|t-s|^{\frac{1}{2} + \frac{3}{2q} }}dt\right)\|\vu(s,\cdot)
\|^2_{L^q}ds\\
=\sup_{\| \psi \|_{L^{p'(\cdot)}}\leq 1}\int_0^T\mathcal{I}_\sigma(|\psi|)(s)\|\vu(s,\cdot)\|^2_{L^q}ds,
\end{eqnarray*}
where $\mathcal{I}_\sigma$ is the 1D Riesz potential given in (\ref{Definition_RieszPotential}) with $\sigma=\frac{1}{2}-\frac{3}{2q}<1$. 
%%%%%%%%%%%%%%%%%%%%%%%%%%%%%%%%%%%%%%%%%%%%%%%%%%%
\begin{Remarque} 
Note that the condition $0<\frac{1}{2}-\frac{3}{2q}<1$ which is needed here in order to use the Riesz potentials implies the constraint $3<q$.
\end{Remarque}
%%%%%%%%%%%%%%%%%%%%%%%%%%%%%%%%%%%%%%%%%%%%%%%%%%%
We apply now the H\"older inequality with $1=\frac{1}{p(\cdot)}+\frac{1}{p(\cdot)}+\frac{1}{ \tilde{p}(\cdot)}$ to obtain the estimate
$$\sup_{\| \psi \|_{L^{p'(\cdot)}}\leq 1}\int_0^T\mathcal{I}_\sigma(|\psi|)(s)\|\vu(s,\cdot)\|^2_{L^q}ds\leq C\sup_{\| \psi \|_{L^{p'(\cdot)}}\leq 1}\|\mathcal{I}_\sigma(|\psi|)\|_{L_t^{\tilde{p}(\cdot)}}\Big\|\|\vu(\cdot,\cdot)\|_{L_x^q}\Big\|_{L_t^{p(\cdot)}}\Big\|\|\vu(\cdot,\cdot)\|_{L_x^q}\Big\|_{L_t^{p(\cdot)}}.$$
%%%%%%%%%%%%%%%%%%%%%%%%%%%%%%%%%%%%%%%%%%%%%%%%%%%
\begin{Remarque} 
The relationship $1=\frac{2}{p(\cdot)}+\frac{1}{ \tilde{p}(\cdot)}$ imposes the condition $p^->2$. 
\end{Remarque}
%%%%%%%%%%%%%%%%%%%%%%%%%%%%%%%%%%%%%%%%%%%%%%%%%%%
By the boundedness of the Riesz potentials in the Lebesgue spaces of variable exponent (see (\ref{PotentialRieszVariable0})) we obtain, 
\begin{eqnarray}
\sup_{\| \psi \|_{L^{p'(\cdot)}}\leq 1}\|\mathcal{I}_\sigma(|\psi|)\|_{L_t^{\tilde{p}(\cdot)}}\|\vu(\cdot,\cdot)\|_{L_t^{p(\cdot)}(L_x^q)}\|\vu(\cdot,\cdot)\|_{L_t^{p(\cdot)}(L_x^q)}\qquad\qquad\qquad\qquad\qquad\notag\\
\leq C\sup_{\| \psi \|_{L^{p'(\cdot)}}\leq 1}\|\psi\|_{L_t^{r(\cdot)}}\|\vu(\cdot,\cdot)\|_{L_t^{p(\cdot)}(L_x^q)}\|\vu(\cdot,\cdot)\|_{L_t^{p(\cdot)}(L_x^q)},\label{Estimation_AvantInclusion}
\end{eqnarray}
where 
\begin{equation}\label{Indices_Riesz}
\frac{1}{\tilde{p}(\cdot)}=\frac{1}{r(\cdot)}-\big(\frac{1}{2}-\frac{3}{2q}\big).
\end{equation}
Note now that if $\frac{2}{p(\cdot)}+\frac{3}{q}<1$ (which is the case by the hypothesis of Theorem \ref{Theoreme_2}), then by the previous identity since we have $\frac{1}{\tilde{p}(\cdot)}=1-\frac{2}{p(\cdot)}$ and $\frac{1}{p'(\cdot)}=1-\frac{1}{p(\cdot)}$, we easily deduce that $r(\cdot)<p'(\cdot)$. This particular relationship will help us to apply the following result
%%%%%%%%%%%%%%%%%%%%%%%%%%%%%%%%%%%%%%%%%%%%%%%%%%%
\begin{Lemme}\label{lemme_embeding}
For $n\geq 1$ and for a given bounded domain $X\subset \mathbb{R}^n$ consider two functions $p_1(\cdot),  p_2(\cdot)\in \mathcal{P}(X)$ such that $1< p_1^+, \ p_2^+ <+\infty$.  Then, $L^{p_2(\cdot)}(X)\subset L^{p_1(\cdot)}(X)$ if and only if $p_1(x)\leq p_2(x)$ almost everywhere. Furthermore, in this case we have that  
$$\|f\|_{L^{p_1(\cdot)} } \leq\left(1+\left|  X \right|\right)\|f\|_{L^{p_2(\cdot)}}.$$
\end{Lemme}
%%%%%%%%%%%%%%%%%%%%%%%%%%%%%%%%%%%%%%%%%%%%%%%%%%%
A proof of this lemma can be found in \cite[Corollary 2.48]{Cruz_Libro}.\\

Thus, if we apply this result to the case $r(\cdot)<p'(\cdot)$ in (\ref{Estimation_AvantInclusion}) we obtain (recall that we are working in the time interval $[0,T]$):
\begin{eqnarray*}
\sup_{\| \psi \|_{L^{p'(\cdot)}}\leq 1}\|\psi\|_{L_t^{r(\cdot)}}\|\vu(\cdot,\cdot)\|_{L_t^{p(\cdot)}(L_x^q)}\|\vu(\cdot,\cdot)\|_{L_t^{p(\cdot)}(L_x^q)}\qquad\qquad\qquad\qquad\qquad\qquad\\
\leq \sup_{\| \psi \|_{L^{p'(\cdot)}}\leq 1}(1+T)\|\psi\|_{L_t^{p'(\cdot)}}\|\vu(\cdot,\cdot)\|_{L_t^{p(\cdot)}(L_x^q)}\|\vu(\cdot,\cdot)\|_{L_t^{p(\cdot)}(L_x^q)}\\
\leq (1+T)\|\vu(\cdot,\cdot)\|_{L_t^{p(\cdot)}(L_x^q)}\|\vu(\cdot,\cdot)\|_{L_t^{p(\cdot)}(L_x^q)}.
\end{eqnarray*}
Now, with all these estimates, getting back to (\ref{Identite_NormeDualite}) we obtain the following control
\begin{equation}\label{Estimation_ApplicationBilineaire}
\left\|\int_{0}^t\mathfrak{g}_{t-s}\ast \mathbb{P}(div(\vu \otimes \vu))(\cdot, \cdot)ds
\right\|_{L_t^{p(\cdot)}(L_x^q)}\leq C(1+T)\|\vu(\cdot,\cdot)\|_{L_t^{p(\cdot)}(L_x^q)}\|\vu(\cdot,\cdot)\|_{L_t^{p(\cdot)}(L_x^q)}.
\end{equation}
\end{itemize}
With the controls (\ref{Estimation_DonneeInitiale}), (\ref{Estimation_ForceExterieure}) and (\ref{Estimation_ApplicationBilineaire}) we have proven the estimates \eqref{Control_Uo_Lplq}, (\ref{Control_force_LpLq}) and \eqref{Control_Nolineal_LpLq}. Now, in order to close the fixed point argument we need the condition
$$C \big(\|\vu_0\|_{L^{q} (\Rt)}+\| \vf\|_{L^1_t(L^{q}_x)}\big)\max\left\{T^{\frac{1}{p^{-}}}, T^{\frac{1}{p^{+}}}\right\}\leq\frac{C}{(1+T)},$$
which can be rewritten in the following manner:
$$\|\vu_0\|_{L^{q} (\Rt)}+\| \vf\|_{L^1_t(L^{q}_x)}\leq\frac{C}{(1+T)\max\left\{T^{\frac{1}{p^{-}}}, T^{\frac{1}{p^{+}}}\right\}},$$
thus, if the initial data and the external force satisfy this estimate, we can obtain a unique (local in time) solution of the Navier-Stokes equations (\ref{NS_Intro}). This concludes the proof of Theorem \ref{Theoreme_2}. \hfill $\blacksquare$

%%%%%%%%%%%%%%%%%%%%%%%%%%%%%%%%%%%%%%%%%%%%%%%%%%%
\begin{Remarque}
Note that if in the estimate \eqref{Estimation_AvantInclusion} we aim to obtain $r(\cdot)=p'(\cdot)$, then by the relationship \eqref{Indices_Riesz} we will obtain the condition $\frac{2}{p(\cdot)} +\frac{3}{q} =1$ that forces the parameter $p(\cdot)$ to be constant.
\end{Remarque}

%%%%%%%%%%%%%%%%%%%%%%%%%%%%%%%%%%%%%%%%%%%%%%%%%%%
\paragraph{\bf Conflict of interest.} On behalf of all authors, the corresponding author states that there is no conflict of interest.\\

\paragraph{\bf Acknowledgements.} 
The second author is supported by the ANID postdoctoral program BCH 2022 grant No. 74220003.
%%%%%%%%%%%%%%%%%%%%%%%%%%%%%%%%%%%%%%%%%%%%%%%%%%%

%%%%%%%%%%%%%%%%%%%%%%%%%

\end{document}